\documentclass[11pt]{article}
\usepackage{amsfonts}
\usepackage{amsmath}
\usepackage{mathrsfs}
\usepackage{color}
\usepackage{mathrsfs,amscd,amssymb,amsthm,amsmath,bm,graphicx,psfrag,subfigure,url}

\makeatletter

\usepackage{indentfirst}

\def \[{\begin{equation}}
\def \]{\end{equation}}

\newtheorem{thm}{Theorem}[section]

\newtheorem{lem}[thm]{Lemma}
\newtheorem{cor}[thm]{Corollary}

\newenvironment{wst}
{\setlength{\leftmargini}{1.5\parindent}
 \begin{itemize}
 \setlength{\itemsep}{-1.1mm}}
{\end{itemize}}

\begin{document}

\title{\bf Relation between the skew-rank of an oriented graph and the independence number of its underlying graph\footnote{S. L. acknowledges the financial support from the National Natural Science Foundation of China (Grant Nos. 11271149, 11371062, 11671164), the Program for New Century Excellent Talents in University (Grant No. NCET-13-0817) and the Special Fund for Basic Scientific Research of Central Colleges (Grant No. CCNU15A02052). H.W. acknowledges the support from Simons Foundation (Grant No. 245307).}}

\author{Jing Huang $^a$ \and Shuchao Li $^{a,}$\footnote{Corresponding author} \and Hua Wang $^b$}

\date{}

\maketitle

\begin{center}
$^a$ Faculty of Mathematics and Statistics, Central China Normal
University, \\ Wuhan 430079, P.R. China\\
1042833291@qq.com (J.~Huang), lscmath@mail.ccnu.edu.cn (S.C.~Li)
\medskip

$^b$ Department of Mathematical Sciences, Georgia Southern University,\\ Statesboro, GA 30460, United States\\
hwang@georgiasouthern.edu (H. ~Wang)
\end{center}


\begin{abstract}
An oriented graph $G^\sigma$ is a digraph without loops or multiple arcs whose underlying graph is $G$. Let $S\left(G^\sigma\right)$ be the skew-adjacency matrix of $G^\sigma$ and $\alpha(G)$ be the independence number of $G$. The rank of $S(G^\sigma)$ is called the skew-rank of $G^\sigma$, denoted by $sr(G^\sigma)$. Wong et al. [European J. Combin. 54 (2016) 76-86] studied the relationship between the skew-rank of an oriented graph and the rank of its underlying graph. In this paper, the correlation involving the skew-rank, the independence number, and some other parameters are considered. First we show that $sr(G^\sigma)+2\alpha(G)\geqslant 2|V_G|-2d(G)$, where $|V_G|$ is the order of $G$ and $d(G)$ is the dimension of cycle space of $G$.  We also obtain sharp lower bounds for $sr(G^\sigma)+\alpha(G),\, sr(G^\sigma)-\alpha(G)$, $sr(G^\sigma)/\alpha(G)$ and characterize all corresponding extremal graphs.
\end{abstract}

\vspace{2mm} \noindent{\bf Keywords}: Skew-rank; Oriented graph; Evenly-oriented; Independence number

\vspace{2mm}


\setcounter{section}{0}
\section{Introduction}\setcounter{equation}{0}

We will start with introducing some background information that will lead to our main results. Some important previously established facts will also be presented.

\subsection{Background}

Let $G=(V_G, E_G)$ be a graph with vertex set $V_G=\{v_1,v_2,\ldots,v_n\}$ and edge set $E_G$. Denote by $P_n, C_n$ and $K_n$ a path, a cycle and a complete graph of order $n$, respectively. The set of neighbors of a vertex $v$ in $G$ is denoted by $N_G(v)$ or simply $N(v).$ Unless otherwise stated, we follow the traditional notations and terminologies (see, for instance, \cite{5}).

The \textit{adjacency matrix} $A(G)$ of $G$ is an $n\times n$ matrix whose $(i, j)$-entry is 1 if vertices $v_i$ and $v_j$ are adjacent and 0 otherwise. Given a graph $G$, the oriented graph $G^\sigma$ is obtained from $G$ by assigning each edge of $G$ a direction. We call $G$ the \textit{underlying graph} of $G^\sigma$. The \textit{skew-adjacency matrix} associated to $G^\sigma$, denoted by $S(G^\sigma)$, is defined to be an $n\times n$ matrix
$[s_{x,y}]$ such that $s_{x,y} = 1$ if there is an arc from $x$ to $y$, $s_{x,y} = -1$ if there is an arc from $y$ to $x$ and $s_{x,y} = 0$
otherwise. The \textit{rank} of $G$, denoted by $r(G)$, is the rank of $A(G)$. The \textit{skew-rank} of $G^\sigma$, denoted by $sr(G^\sigma)$, is the rank of $S(G^\sigma)$. It is easy to see that $sr(G^\sigma)$ is even since $S(G^\sigma)$ is skew symmetric.

The value
$$d(G):=|E_G|-|V_G|+\omega(G),$$
is called the \textit{dimension} of cycle space of $G$, where $\omega(G)$ is the number of the components of $G$. Two distinct edges in a graph $G$ are \textit{independent} if they do not share a common end-vertex. A \textit{matching} is a set of pairwise independent edges of $G$, while a \textit{maximum matching} of $G$ is a matching with the maximum cardinality. The \textit{matching number} of $G$, written as $\alpha'(G)$, is the cardinality of a maximum matching of $G.$ Two vertices of a graph $G$ are said to be \textit{independent} if they are not adjacent. A subset $I$ of $V_G$ is called an \textit{independent set} if any two vertices of $I$ are independent in $G$. An independent set $I$ is \textit{maximum} if $G$ has no independent set $I'$ with $|I'|>|I|$. The number of vertices in a maximum independent set of $G$ is called the \textit{independence number} of $G$ and is denoted by $\alpha(G)$. An oriented graph is called  \textit{acyclic} (resp. \textit{connected,\, bipartite}) if its underlying graph is acyclic (resp. connected,\, bipartite). A graph is called an \textit{empty graph} if it has no edges. We call $v$ a \textit{cut-vertex} of a connected $G^\sigma$ if $G^\sigma-v$ is disconnected.

The study on skew spectrum of oriented graphs has attracted much attention. Anuradha and Balakrishnan \cite{02} investigated skew spectrum of the Cartesian product of two oriented graphs. Anuradha et al. \cite{03} considered the skew spectrum of bipartite graphs. Hou and Lei \cite{015} studied the coefficients of the characteristic polynomial of skew-adjacency matrix of a oriented graph. 
Xu \cite{025} established a relation between the spectral radius and the skew spectral radius. Cavers et al. \cite{05} systematically studied skew-adjacency matrices of directed graphs. Among various specific topics, the minimal skew-rank of oriented graphs is of particular interest to researchers. The graphs with minimum skew rank 2, 4 are characterized in \cite{006}. The bicyclic oriented graphs with skew-rank 2, 4 and 6 are, respectively, determined in \cite{8} and \cite{9}. Recently, Mallik and Shader \cite{10} studied the minimum rank of all real skew-symmetric matrices described by a graph. Wong, Ma and Tian \cite{1} presented a beautiful relation between the skew-rank of an oriented graph and the rank of its underlying graph. Huang and Li \cite{004} further extended these results. For more properties and applications of the skew-rank of oriented graphs, we refer the readers to \cite{12,13,14,11}.

Very recently, Ma, Wong and Tian \cite{6} determined the relationship between $sr(G^\sigma)$ and the matching number $\alpha'(G)$. They \cite{31} also characterized the relationship between $r(G)$ and $p(G)$ (the number of pendant vertices of $G$), from which one can obtain the same relationship between $sr(G^\sigma)$ and $p(G)$. It is natural to further this study by considering the correlation between $sr(G^\sigma)$ and some other parameters of its underlying graph. In this paper we first establish the sharp lower bound on $sr(G^\sigma)+2\alpha(G)$ of an oriented graph. We then apply the same fundamental idea to determine sharp lower bounds on $sr(G^\sigma)+\alpha(G), sr(G^\sigma)-\alpha(G)$ and $sr(G^\sigma)/\alpha(G)$ and characterize the corresponding extremal oriented graphs.

\subsection{Main results}

Let $C_k=v_1v_2\cdots v_kv_1$ be a cycle of length $k$. The \textit{sign} of $C_k^\sigma$ with respect to $\sigma$ is defined to be the
sign of $\left(\prod_{i=1}^{k-1}s_{i,i+1}\right)\cdot s_{k,1}.$ An even oriented cycle $C_k^\sigma$ is called \textit{evenly-oriented} (resp. \textit{oddly-oriented}) if its sign is positive (resp. negative). An induced subgraph of $G^\sigma$ is an induced subgraph of $G$ where each edge preserves the original orientation in $G^\sigma$. For an induced subgraph $H^\sigma$ of $G^\sigma$, let $G^\sigma-H^\sigma$ be the subgraph obtained from $G^\sigma$ by removing all vertices of $H^\sigma$ and their incident edges. For $W\subseteq V_{G^\sigma}$, $G^\sigma-W$ is the subgraph obtained from $G^\sigma$ by removing all vertices in $W$ and all incident edges. A vertex of $G^\sigma$ is called a \textit{pendant vertex} if it is of degree one in $G$, whereas a vertex of $G^\sigma$ is called a \textit{quasi-pendant vertex} if it is adjacent to a pendant vertex in $G$.

Given a graph $G$ with pairwise vertex-disjoint cycles, let $\mathscr{C}_G$ denote the set of all cycles of $G$. Contracting each cycle to a single vertex yields an acyclic graph $T_G$ from $G$.
It is clear that $T_G$ is always acyclic. Note that the graph $T_G-W_{\mathscr{C}}$ (where $W_{\mathscr{C}}$ is the set of vertices corresponding to the cycles in $G$)  is the same as the graph obtained from $G$ by removing all the vertices on cycles and their incident edges. We denote this graph by $\Gamma_G.$ For example, in Fig. 1, $T_G$ is obtained from $G$ by contracting each cycle into a single vertex, and $\Gamma_G$ is obtained from $G$ by removing all the vertices on cycles and their incident edges.
\begin{figure}[h!]
\begin{center}
\psfrag{a}{$G$}\psfrag{b}{$T_G$}
\psfrag{c}{$\Gamma_G$}
\includegraphics[width=120mm]{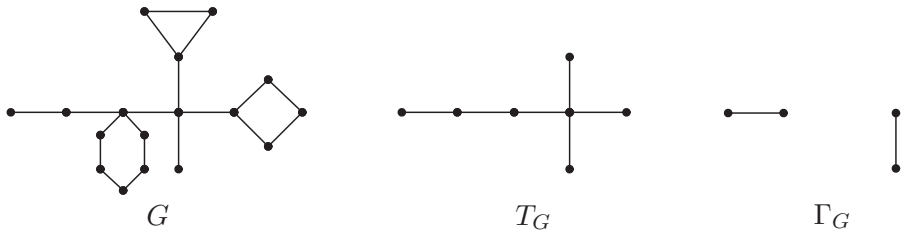} \\
  \caption{ Graphs $G,\, T_G$, and $\Gamma_G$.} \label{fig:1}
\end{center}
\end{figure}

Following the above notations our first main result reads as follows. 
\begin{thm}\label{theo:2.1}
Let $G^\sigma$ be a simple connected graph on $n$ vertices. Then
\begin{equation}\label{eq:1.1}
sr(G^\sigma)+2\alpha(G)\geqslant2n-2d(G).
\end{equation}
The equality in $(\ref{eq:1.1})$ holds if and only if the following conditions hold for $G^\sigma:$
\begin{wst}
\item[{\rm (i)}] the cycles (if any) of $G^\sigma$ are pairwise vertex-disjoint;
\item[{\rm (ii)}] each cycle of $G^\sigma$ is odd or evenly-oriented;
\item[{\rm (iii)}]$\alpha(T_G)=\alpha(\Gamma_G)+d(G)$.
\end{wst}
\end{thm}
For example, let $G$ be as in Fig.~\ref{fig:1}. If all the even cycles in $G^\sigma$ are evenly-oriented, then $G^\sigma$ satisfies conditions (i)-(iii) (note that $\alpha(T_G)=5,\, \alpha(\Gamma_G)=2,\,d(G)=3)$ and $sr(G^\sigma)+2\alpha(G)=2n-2d(G)$ holds with $n=17,\,sr(G^\sigma)=12$ and $\alpha(G)=8$.

In the case that $G$ is bipartite, the following is a direct consequence of Theorem~\ref{theo:2.1}.
\begin{cor}\label{cor:2.2}
Let $G^\sigma$ be a simple connected bipartite graph with $n$ vertices. Then $sr(G^\sigma)+2\alpha(G)\geqslant2n-2d(G)$ with equality if and only if the following conditions hold for $G^\sigma:$
\begin{wst}
\item[{\rm (i)}] the cycles (if any) of $G^\sigma$ are pairwise vertex-disjoint;
\item[{\rm (ii)}] each cycle of $G^\sigma$ is evenly-oriented;
\item[{\rm (iii)}]$\alpha(T_G)=\alpha(\Gamma_G)+d(G)$.
\end{wst}
\end{cor}

Note that $\alpha(G)+\alpha'(G)=|V_G|$ if $G$ is bipartite and $d(G)$ is exactly the number of cycles if the cycles of $G$ are pairwise vertex-disjoint. Then Corollary~\ref{cor:2.2} is equivalent to Theorem~\ref{theo:2.3} below when $G$ is bipartite, obtained in
\cite{6}, showing the correlation between the skew-rank of an oriented graph, the matching number, and the dimension of cycle space of its underlying graph.

\begin{thm}[\cite{6}]\label{theo:2.3}
Let $G^\sigma$ be a simple connected graph. Then $sr(G^\sigma)-2\alpha'(G)\geqslant-2d(G)$ with equality if and only if the following conditions hold for $G^\sigma:$
\begin{wst}
\item[{\rm (i)}] the cycles (if any) of $G^\sigma$ are pairwise vertex-disjoint;
\item[{\rm (ii)}] each cycle of $G^\sigma$ is evenly-oriented;
\item[{\rm (iii)}] $\alpha'(T_G)=\alpha'(\Gamma_G)$.
\end{wst}
\end{thm}

Along the same line, we establish sharp lower bounds on $sr(G^\sigma)+\alpha(G), sr(G^\sigma)-\alpha(G)$, and $sr(G^\sigma)/\alpha(G)$ in the next three theorems.

\begin{thm}\label{theo:2.4}
Let $G^\sigma$ be a simple connected graph with $n$ vertices and $m$ edges. Then
\begin{equation}\label{eq:1.2}
sr(G^\sigma)+\alpha(G)\geqslant 4n-2m-\sqrt{n(n-1)-2m+\frac{1}{4}}-\frac{5}{2}
\end{equation}
with equality if and only if $G\cong S_n$ or $G\cong C_3$.
\end{thm}

\begin{thm}\label{theo:2.5}
Let $G^\sigma$ be a simple connected graph with $n$ vertices and $m$ edges. Then
\begin{equation*}
sr(G^\sigma)-\alpha(G)\geqslant 4n-2m-3\sqrt{n(n-1)-2m+\frac{1}{4}}-\frac{7}{2}
\end{equation*}
with equality if and only if $G\cong S_n$ or $G\cong C_3$.
\end{thm}

\begin{thm}\label{theo:2.6}
Let $G^\sigma$ be a simple connected graph with $n$ vertices and $m$ edges. Then
\begin{equation*}
\frac{sr(G^\sigma)}{\alpha(G)}\geqslant \frac{4(2n-m-1)}{\sqrt{4n(n-1)-8m+1}+1}-2
\end{equation*}
with equality if and only if $G\cong S_n$ or $G\cong C_3$.
\end{thm}

In the rest of this section we recall some important known results. In Section~\ref{sec:2} we first establish some technical lemmas that help us characterize the extremal graphs. We present the proofs of our main results in Section~\ref{sec:3}. We briefly comment on our findings and propose some questions in Section~\ref{sec:4}.

\subsection{Preliminaries}

For the rest of our introduction we recall the following important facts.

\begin{lem}[\cite{2}]
\label{lem:3.1}

Let $G^\sigma$ be an oriented graph:
\begin{wst}
\item[{\rm (i)}] If $H^\sigma$ is an induced subgraph of $G^\sigma$, then $sr(H^\sigma)\leqslant sr(G^\sigma);$
\item[{\rm (ii)}] If $G_1^\sigma, G_2^\sigma, \ldots, G_t^\sigma$ are all the components of $G^\sigma,$ then $sr(G^\sigma)=\sum^t_{i=1}sr(G_i^\sigma);$
\item[{\rm (iii)}] $sr(G^\sigma)\geqslant 0$ with equality if and only if $G^\sigma$ is an empty graph.
\end{wst}
\end{lem}

The following observation immediately follows from the definition of the independence number.

\begin{lem}\label{lem:3.2}
Let $G$ be a simple connected graph. Then
\begin{wst}
\item[{\rm (i)}]$\alpha(G)-1\leqslant\alpha(G-v)\leqslant\alpha(G)$ for any $v\in V_G$;
\item[{\rm (ii)}] $\alpha(G-e)\geqslant\alpha(G)$ for any $e\in E_G$.
\end{wst}
\end{lem}

\begin{lem}[\cite{5}]
\label{lem:3.3}
Let $P_n$ be a path of order $n$. Then $r(P_n)=n$ if $n$ is even, and $r(P_n)=n-1$ if $n$ is odd.
\end{lem}

\begin{lem}[\cite{2}]
\label{lem:3.4}
Let $F^\sigma$ be an oriented acyclic graph with matching number $\alpha'(F)$. Then $sr(F^\sigma)=r(F)=2\alpha'(F).$
\end{lem}

\begin{lem}[\cite{28}]
\label{lem:3.5}
Let $G$ be a bipartite graph with $n$ vertices. Then $\alpha(G)+\alpha'(G)=n$.
\end{lem}

\begin{lem}[\cite{23}]
\label{lem:3.6}
Let $C_n^\sigma$ be an oriented cycle of order $n$. Then $sr(C_n^\sigma)=n$ if $C_n^\sigma$ is oddly-oriented, $sr(C_n^\sigma)=n-2$ if $C_n^\sigma$ is evenly-oriented and $sr(C_n^\sigma)=n-1$ if $n$ is odd.
\end{lem}

\begin{lem}[\cite{2}]
\label{lem:3.7}
Let $y$ be a pendant vertex of $G^\sigma,$\ and $x$ be the neighbor of $y$, then $sr(G^\sigma)=sr(G^\sigma-x)+2=
sr(G^\sigma-x-y)+2$.
\end{lem}

\begin{lem}[\cite{1}]
\label{lem:3.8}
Let $x$ be a vertex of $G^\sigma$. Then $sr(G^\sigma-x)$ is equal to either $sr(G^\sigma)$ or $sr(G^\sigma)-2.$
\end{lem}

The following lemma on the dimension of cycle space of $G$ follows directly from the definition of $d(G)$.

\begin{lem}[\cite{1}]
\label{lem:3.9}

Let $G$ be a graph with $x\in V_G$.
\begin{wst}
\item[{\rm (i)}] $d(G)=d(G-x)$ if $x$ is not on any cycle of $G;$
\item[{\rm (ii)}] $d(G-x)\leqslant d(G)-1$ if $x$ lies on a cycle;
\item[{\rm (iii)}] $d(G-x)\leqslant d(G)-2$ if $x$ is a common vertex of distinct cycles;
\item[{\rm (iv)}] If the cycles of $G$ are pairwise vertex-disjoint, then $d(G)$ is exactly the number of cycles in $G$.
\end{wst}
\end{lem}

The next result is on the rank of an acyclic graph. Let $T$ be a tree with at least one edge, we denote by $\widetilde{T}$ the subtree obtained from $T$ by removing all pendant vertices of $T$.

\begin{lem}[\cite{6}]
\label{lem:3.10}
Let $T$ be a tree with at least one edge. Then
\begin{wst}
\item[{\rm (i)}] $r(\widetilde{T})<r(T);$
\item[{\rm (ii)}] If $r(T-D)=r(T)$ for a subset $D$ of $V_T$, then there is a pendant vertex $v$ such that $v\notin D.$
\end{wst}
\end{lem}

Recall that $p(G)$ is the number of pendant vertices of $G$, from Lemmas~\ref{lem:3.4}, \ref{lem:3.5} and \ref{lem:3.10} we immediately have the following.

\begin{cor}
\label{cor:3.11}
Let $T$ be a tree with at least one edge. Then
\begin{wst}
\item[{\rm (i)}] $\alpha(T)<\alpha(\widetilde{T})+p(T)$;
\item[{\rm (ii)}] If $\alpha(T)=\alpha(T-D)+|D|$ for a subset $D$ of $V_T$, then there is a pendant vertex $v$ such that $v\notin D.$
\end{wst}
\end{cor}

\section{Technical lemmas}\label{sec:2}

In this section we present a few technical lemmas. First we establish \eqref{eq:1.1}.

\begin{lem}\label{lem:new}
The inequality \eqref{eq:1.1} holds.
\end{lem}

\begin{proof}
We proceed by induction on $d(G).$  If $d(G)=0$, then $G^\sigma$ is an oriented tree and the result follows immediately from Lemmas~\ref{lem:3.4} and \ref{lem:3.5}. Now suppose that $G^\sigma$ has at least one cycle, i.e., $d(G) \geqslant 1,$ and let $x$ be a vertex on some cycle.  By Lemma~\ref{lem:3.9}(ii) we have
\begin{eqnarray}\label{eq:3.1}
d(G-x)\leqslant d(G)-1.
\end{eqnarray}
By the induction hypothesis one has
\begin{eqnarray}\label{eq:3.2}
sr(G^\sigma-x)+2\alpha(G-x)\geqslant2(n-1)-2d(G-x).
\end{eqnarray}
By Lemma~\ref{lem:3.1}(i) and Lemma~\ref{lem:3.2}(i), we obtain
\begin{eqnarray}\label{eq:3.3}
 sr(G^\sigma-x)\leqslant sr(G^\sigma),\ \ \ \alpha(G-x)\leqslant\alpha(G).
\end{eqnarray}
The inequality \eqref{eq:1.1} then follows from (\ref{eq:3.1})-(\ref{eq:3.3}).
\end{proof}

For convenience we call a graph $G^\sigma$ ``lower optimal'' if it achieves equality in \eqref{eq:1.1}. In the rest of this section we aim to provide some fundamental characterizations of lower-optimal oriented graphs.

\begin{lem}\label{lem:4.1}
Let $x$ be a vertex on a cycle of $G^\sigma$. If $G^\sigma$ is lower-optimal, then
\begin{wst}
\item[{\rm (i)}] $sr(G^\sigma)=sr(G^\sigma-x);$
\item[{\rm (ii)}] $\alpha(G)=\alpha(G-x);$
\item[{\rm (iii)}] $d(G)=d(G-x)+1;$
\item[{\rm (iv)}] $G^\sigma-x$ is lower-optimal;
\item[{\rm (v)}] $x$ lies on just one cycle of $G$ and $x$ is not a quasi-pendant vertex of $G$.
\end{wst}
\end{lem}

\begin{proof}
The lower-optimal condition for $G^\sigma$ together with the proof of Lemma~\ref{lem:new} forces equalities in (\ref{eq:3.1})-(\ref{eq:3.3}). Consequently we have (i)-(iv). By (iii) and Lemma~\ref{lem:3.9}(iii) we obtain that $x$ lies on just one cycle of $G$. If $x$ is a quasi-pendant vertex adjacent to a pendant vertex $y$, then by Lemma~\ref{lem:3.7}, we have $sr(G^\sigma)=sr(G^\sigma-x)+2$, which is a contradiction to (i). This completes the proof of (v).
\end{proof}

The next observation, although simple, is very helpful to our proof.

\begin{lem}\label{lem:4.2}
Let $y$ be a pendant vertex of $G$ with neighbor $x$. Then $\alpha(G)=\alpha(G-x)=\alpha(G-x-y)+1$.
\end{lem}

\begin{proof}
It is routine to check that $\alpha(G-x)=\alpha(G-x-y)+1$. In order to complete the proof, it suffices to show that $\alpha(G)=\alpha(G-x)$. In fact, let $I$ be a maximum independent set of $G$.

If $x\notin I$, then $I$ is also a maximum independent set of $G-x$ and we have $\alpha(G)=|I|=\alpha(G-x)$.

If $x\in I$, then $y\notin I$, thus $(I\backslash \{x\})\cup\{y\}$ is an independent set of $G-x$. Hence we have $\alpha(G-x)\geqslant\left|(I\backslash \{x\})\cup\{y\}\right|=|I|=\alpha(G)$. By Lemma~\ref{lem:3.2}(i), we have $\alpha(G-x)\leqslant\alpha(G)$.

Therefore we have $\alpha(G)=\alpha(G-x)=\alpha(G-x-y)+1$ as desired.
\end{proof}

Given an induced oriented subgraph $H^\sigma$ of $G^\sigma,$  let $v_i$ be in $V_{G^\sigma} \setminus V_{H^\sigma}$. Then the induced oriented subgraph of $G^\sigma$ with vertex set $V_{H^\sigma}\bigcup\{v_i\}$ is simply written as $H^\sigma+v_i$. The following lemma summarizes a few known results.

\begin{lem}\label{lem:4.3}
Let $C_q^\sigma$ be a pendant oriented cycle of $G^\sigma$ with $x$ being a vertex of $C_q$ of degree $3,$ and let $H^\sigma =
G^\sigma-C_q^\sigma, M^\sigma = H^\sigma + x.$ Then
\begin{eqnarray*}
sr(G^\sigma)=\left\{
                             \begin{array}{lll}
                               q-1+sr(M^\sigma), & \hbox{if $q$ is odd;}  \ \ \ \ \ \text{\rm (see \cite{1})}\\
                               q-2+sr(M^\sigma), & \hbox{if $C_q^\sigma$ is evenly-oriented;} \ \ \ \ \ \text{\rm (see \cite{1})} \\
                               q+sr(H^\sigma), & \hbox{if $C_q^\sigma$ is oddly-oriented.} \ \ \ \ \ \text{\rm (see \cite{004})}
                               \end{array}
                            \right.
\end{eqnarray*}
\end{lem}

Following the same direction we establish a few more facts in the rest of this section.

\begin{lem}\label{lem:4.4}
Let $C_q^\sigma$ be a pendant oriented cycle of $G^\sigma$ with $x$ being the unique vertex of $C_q$ of degree 3. Let $H^\sigma=
G^\sigma-C_q^\sigma$ and  $M^\sigma = H^\sigma + x$, if $G^\sigma$ is lower-optimal, then
\begin{wst}
\item[{\rm (i)}] $q$ is odd or $C_q^\sigma$ is evenly-oriented;
\item[{\rm (ii)}] $sr(G^\sigma)=q-1+sr(H^\sigma), \alpha(G)=\alpha(H)+\frac{q-1}{2}$ if $q$ is odd and $sr(G^\sigma)=q-2+sr(H^\sigma), \alpha(G)=\alpha(H)+\frac{q}{2}$ if $C_q^\sigma$ is evenly-oriented;
\item[{\rm (iii)}] both $H^\sigma$ and $M^\sigma$ are lower-optimal;
\item[{\rm (iv)}] $sr(M^\sigma)=sr(H^\sigma)$ and $\alpha(M) = \alpha(H)+1.$
\end{wst}
\end{lem}

\begin{proof}
(i)\ Supposing for contradiction that $C_q^\sigma$ is oddly-oriented, then by Lemma~\ref{lem:4.3} we have
\begin{eqnarray}\label{eq:3.4}
sr(G^\sigma)=q+sr(H^\sigma).
\end{eqnarray}
Note that $x$ lies on the cycle $C_q$. Hence, by Lemma~\ref{lem:4.1}(ii) we have
\begin{eqnarray}\label{eq:3.5}
\alpha(G)=\alpha(G-x)=\alpha(P_{q-1})+\alpha(H)=\frac{q}{2}+\alpha(H).
\end{eqnarray}
As $C_q$ is a pendant cycle of $G$, we have
\begin{eqnarray}\label{eq:3.6}
d(G)=d(M)+1=d(H)+1.
\end{eqnarray}
Suppose $|V_G|=n$. Since $G^\sigma$ is lower-optimal, we have
\begin{eqnarray}\label{eq:3.7}
sr(G^\sigma)+2\alpha(G)=2n-2d(G).
\end{eqnarray}
From (\ref{eq:3.4})-(\ref{eq:3.7}) we have $sr(H^\sigma)+2\alpha(H)=2(n-q)-2d(H)-2$, which is a contradiction to (\ref{eq:1.1}).
This completes the proof of (i).

Next we show (ii)-(iv) according to the following two possible cases.\vspace{3mm}

\noindent {\bf{Case 1.}}\ $q$ is odd.\vspace{3mm}

(ii)\ Note that $x$ lies on a cycle of $G$, by Lemma~\ref{lem:4.1}(i)-(ii) we have
\begin{align}
sr(G^\sigma)&=sr(G^\sigma-x)=sr(P_{q-1}^\sigma)+sr(H^\sigma)=q-1+sr(H^\sigma),\label{eq:3.8}\\
\alpha(G)&=\alpha(G-x)=\alpha(P_{q-1})+\alpha(H)=\frac{q-1}{2}+\alpha(H).\label{eq:3.9}
\end{align}

(iii)-(iv)\ From (\ref{eq:3.6})-(\ref{eq:3.9}) we have $sr(H^\sigma)+2\alpha(H)=2(n-q)-2d(H)$, implying that $H^\sigma$ is lower-optimal.

Since $q$ is odd, by Lemma~\ref{lem:4.3}, we have
\begin{eqnarray}\label{eq:3.10}
sr(G^\sigma)=q-1+sr(M^\sigma).
\end{eqnarray}
Combining (\ref{eq:3.8}) and (\ref{eq:3.10}) yields
\begin{eqnarray}\label{eq:3.11}
sr(H^\sigma)=sr(M^\sigma).
\end{eqnarray}

Furthermore,
\begin{eqnarray}
2\alpha(H)&=&2(n-q)-sr(H^\sigma)-2d(H)\notag\\
&=&2(n-q+1)-sr(M^\sigma)-2d(H)-2\notag\\
&=&2(n-q+1)-sr(M^\sigma)-2d(M)-2\notag\\
&\leqslant&2\alpha(M)-2,\label{eq:3.12}
\end{eqnarray}
where the first equality follows from the lower-optimal condition for $H^\sigma$, the second and the third equalities follow from (\ref{eq:3.11}) and (\ref{eq:3.6}), respectively. And the last inequality (\ref{eq:3.12}) follows from applying (\ref{eq:1.1}) to $M^\sigma$. Thus we have $\alpha(H)\leqslant\alpha(M)-1$. It follows from Lemma~\ref{lem:3.2}(i) that $\alpha(H)\geqslant\alpha(M)-1$. Hence $\alpha(H)=\alpha(M)-1$.

Consequently we have $sr(M^\sigma)+2\alpha(M)=2(n-q+1)-2d(M)$, implying that $M^\sigma$ is also lower-optimal.

\noindent {\bf{Case 2.}}\ $C_q^\sigma$ is evenly-oriented.

(ii)\ Since $x$ lies on a cycle of $G$, by Lemma~\ref{lem:4.1}(i)-(ii) we have
\begin{align}
sr(G^\sigma)&=sr(G^\sigma-x)=sr(P_{q-1}^\sigma)+sr(H^\sigma)=q-2+sr(H^\sigma),\label{eq:3.13}\\
\alpha(G)&=\alpha(G-x)=\alpha(P_{q-1})+\alpha(H)=\frac{q}{2}+\alpha(H).\label{eq:3.14}
\end{align}

(iii)\ Let $x_1$ be on $C_q$ such that it is adjacent to $x$. By applying Lemma~\ref{lem:4.1} to $G^\sigma$ (resp.  $G$) and Lemma~\ref{lem:3.7} (resp. Lemma~\ref{lem:4.2}) to $G^\sigma-x_1$ (resp.  $G-x_1$) we have
\begin{align}
sr(G^\sigma)&=sr(G^\sigma-x_1)=q-2+sr(M^\sigma),\label{eq:3.15}\\
\alpha(G)&=\alpha(G-x_1)=\frac{q-2}{2}+\alpha(M)\label{eq:3.16}.
\end{align}

From (\ref{eq:3.6})-(\ref{eq:3.7}) and (\ref{eq:3.13})-(\ref{eq:3.14}), one has $sr(H^\sigma)+2\alpha(H)=2(n-q)-2d(H)$, implying that $H^\sigma$ is lower-optimal.

Combining  (\ref{eq:3.6})-(\ref{eq:3.7}) and (\ref{eq:3.15})-(\ref{eq:3.16}), we have $sr(M^\sigma)+2\alpha(M)=2(n-q+1)-2d(M)$, which implies that $M^\sigma$ is also lower-optimal.

(iv)\ \ Combining (\ref{eq:3.13}) and (\ref{eq:3.15}) yields $sr(M^\sigma)=sr(H^\sigma),$ whereas equalities (\ref{eq:3.14}) and (\ref{eq:3.16}) lead to $\alpha(M)=\alpha(H)+1.$

This completes the proof.
\end{proof}

\begin{lem}\label{lem:4.5}
Let $y$ be a pendant vertex of $G^\sigma$ with neighbor $x$, and let $H^\sigma=G^\sigma-y-x.$ If $G^\sigma$ is lower-optimal, then
\begin{wst}
\item[{\rm (i)}] $x$ does not lie on any cycle of $G;$
\item[{\rm (ii)}] $H^\sigma$ is also lower-optimal.
\end{wst}
\end{lem}
\begin{proof}
(i)\ \ Since $x$ is a quasi-pendant vertex of $G$, Lemma~\ref{lem:4.1}(v) states that $x$ does not lie on any cycle of $G.$

(ii)\ \ By Lemmas~\ref{lem:3.7} and \ref{lem:4.2}, we have
\begin{eqnarray}
sr(G^\sigma)=sr(H^\sigma)+2 \hbox{ and } \alpha(G)=\alpha(H)+1.\label{eq:4.17}
\end{eqnarray}
Since $x$ does not lie on any cycle of $G$, by Lemma~\ref{lem:3.9}(i) we have
\begin{eqnarray}\label{eq:4.18}
d(G)=d(H).
\end{eqnarray}
Equalities (\ref{eq:4.17})-(\ref{eq:4.18}), together with the lower-optimal condition of $G^\sigma$, imply that $sr(H^\sigma)+2\alpha(H)=2(n-2)-2d(H)$, i.e., $H^\sigma$ is lower-optimal.
\end{proof}

\begin{lem}\label{lem:4.6}
If $G^\sigma$ is lower-optimal, then
\begin{wst}
\item[{\rm (i)}] the cycles (if any) of $G^\sigma$ are pairwise vertex-disjoint;
\item[{\rm (ii)}] each cycle (if any) of $G^\sigma$ is odd or evenly-oriented;
\item[{\rm (iii)}] $\alpha(G)=\alpha(T_G)+\sum_{C\in \mathscr{C}_G}\left\lfloor\frac{|V_C|}{2}\right\rfloor-d(G).$
\end{wst}
\end{lem}
\begin{proof}
If $G$ contains cycles, then let $x$ be a vertex on some cycle. By Lemma~\ref{lem:4.1}(iii) we have $d(G)=d(G-x)+1$. By Lemma~\ref{lem:3.9}(iii) $x$ can not be a common vertex of distinct cycles, hence the cycles of $G^\sigma$ are pairwise vertex-disjoint. This completes the proof of (i).

We proceed by induction on the order $n$ of $G$ to prove (ii) and (iii). The initial case $n=1$ is trivial. Suppose that (ii) and (iii) hold for any lower-optimal oriented graph of order smaller than $n$, and suppose $G^\sigma$ is an lower-optimal oriented graph of order $n\geqslant 2.$

If $T_G$ is an empty graph, then $G^\sigma$ is a single oriented cycle. Thus (ii) follows from the fact that a single oriented cycle $C_q^\sigma$ is lower-optimal if and only if $q$ is odd or $C_q^\sigma$ itself is evenly-oriented. And (iii) follows from the fact that $\alpha(C_q)=\frac{q-1}{2}$ if $q$ is odd and $\alpha(C_q)=\frac{q}{2}$ if $C_q^\sigma$ is evenly-oriented.

If $T_G$ has at least one edge, then $T_G$ contains at least one pendant vertex, say $y$. Then $y$ is either a pendant vertex of $G$ or  $y\in W_{\mathscr{C}}$, in which case $G$ contains a pendant cycle. We now consider both cases.\vspace{3mm}

\noindent {\bf{Case 1.}}\ $G$ contains a pendant vertex $y$. In this case, let $x$ be the neighbor of $y$ in $G$ and let $H^\sigma=G^\sigma-x-y.$ By Lemma~\ref{lem:4.5}, $x$ is not a vertex on any cycle of $G$ and $H^\sigma$ is also lower-optimal. By induction hypothesis we have
\begin{wst}
\item[{\rm (a)}] each cycle (if any) of $H^\sigma$ is odd or evenly-oriented;
\item[{\rm (b)}] $\alpha(H)=\alpha(T_H)+\sum_{C\in \mathscr{C}_H}\left\lfloor\frac{|V_C|}{2}\right\rfloor-d(H).$
\end{wst}

Note that all cycles of $G$ are also in $H$, Assertion (a) implies that each cycle (if any) of $G^\sigma$ is odd or evenly-oriented. Hence (ii) holds.

Since $x$ does not lie on any cycle of $G$, by Lemma~\ref{lem:3.9}(i) we have
\begin{eqnarray}\label{eq:3.17}
d(G)=d(H).
\end{eqnarray}
Recall that $y$ is also a pendant vertex of $T_G$ adjacent to $x$ and $T_H=T_G-x-y$, then by Lemma~\ref{lem:4.2}, Assertion (b) and (\ref{eq:3.17}) we have
$$
\alpha(G)=\alpha(H)+1=\alpha(T_H)+\sum_{C\in \mathscr{C}_H}\left\lfloor\frac{|V_C|}{2}\right\rfloor-d(H)+1=\alpha(T_G)+\sum_{C\in \mathscr{C}_G}\left\lfloor\frac{|V_C|}{2}\right\rfloor-d(G).
$$
Thus (iii) holds.\vspace{3mm}

\noindent {\bf{Case 2.}}\ $G$ has a pendant cycle $C_q$. In this case, let $x$ be the unique vertex of $C_q$ of degree 3, $H^\sigma=G^\sigma-C_q^\sigma$ and $M^\sigma=H^\sigma+x.$ It follows from Lemma~\ref{lem:4.4}(iii) that $M^\sigma$ is lower-optimal. Applying the induction hypothesis to $M^\sigma$ yields
\begin{wst}
\item[{\rm (c)}] each cycle of $M^\sigma$ is odd or evenly-oriented;
\item[{\rm (d)}] $\alpha(M)=\alpha(T_M)+\sum_{C\in \mathscr{C}_M}\left\lfloor\frac{|V_C|}{2}\right\rfloor-d(M).$
\end{wst}
Assertion (c) and Lemma~\ref{lem:4.4}(i) imply that each cycle of $G^\sigma$ is odd or evenly-oriented since $\mathscr{C}_G=\mathscr{C}_M\bigcup{C_q}.$ Thus, (ii) holds.

Combining Lemma~\ref{lem:4.4}(ii), Lemma~\ref{lem:4.4}(iv) and Assertion (d) we have
\begin{eqnarray}\label{eq:3.18}
\alpha(G)=\alpha(M)+\left\lfloor\frac{q}{2}\right\rfloor-1=\alpha(T_M)+\sum_{C\in \mathscr{C}_M}\left\lfloor\frac{|V_C|}{2}\right\rfloor+\left\lfloor\frac{q}{2}\right\rfloor-d(M)-1.
\end{eqnarray}
As $C_q$ is a pendant cycle of $G$, we have
\begin{eqnarray}\label{eq:3.19}
d(G)=d(M)+1.
\end{eqnarray}
Note that $T_M\cong T_G$ and $\left\lfloor\frac{q}{2}\right\rfloor+\sum_{C\in \mathscr{C}_M}\left\lfloor\frac{|V_C|}{2}\right\rfloor=\sum_{C\in \mathscr{C}_G}\left\lfloor\frac{|V_C|}{2}\right\rfloor$. Together with (\ref{eq:3.18})-(\ref{eq:3.19}) we have
\begin{eqnarray*}
\alpha(G)=\alpha(T_G)+\sum_{C\in \mathscr{C}_G}\left\lfloor\frac{|V_C|}{2}\right\rfloor-d(G),
\end{eqnarray*}
as desired.
\end{proof}

\section{Proofs of main results}\label{sec:3}\setcounter{equation}{0}

We will first provide the proof of Theorem~\ref{theo:2.1}, based on which the other proofs follow.

\subsection{Theorem~\ref{theo:2.1}}

Lemma~\ref{lem:new} already established (\ref{eq:1.1}). We now characterize all the oriented graphs $G^\sigma$ which attain the lower bound by considering the sufficient and necessary conditions for the equality in (\ref{eq:1.1}).

For ``sufficiency", we proceed by induction on the order $n$ of $G$ to show that $G^\sigma$ is lower-optimal if $G^\sigma$ satisfies the conditions (i)-(iii). The $n=1$ case is trivial. Suppose that any oriented graph of order smaller than $n$ which satisfies (i)-(iii) is lower-optimal, and suppose $G^\sigma$ is an oriented graph with order $n\geqslant 2$ that satisfies (i)-(iii). Since the cycles (if any) of $G^\sigma$ are pairwise vertex-disjoint, Lemma~\ref{lem:3.9}(iv) states that $G$ has exactly $d(G)$ cycles, implying that $|W_{\mathscr{C}}|=d(G)$.  If $T_G$ is an empty graph, it follows from (ii) that $G^\sigma$ is an odd cycle or an evenly-oriented cycle, leading to the fact that $G^\sigma$ is lower-optimal. So in what follows, we assume that $T_G$ has at least one edge.

Note that $\alpha(T_G)=\alpha(\Gamma_G)+d(G)=\alpha(T_G-W_{\mathscr{C}})+d(G)$. Then by Corollary~\ref{cor:3.11}(ii), there exists a pendant vertex of $T_G$ not in $W_{\mathscr{C}}$. Thus $G$ contains as least one pendant vertex, say $y$. Let $x$ be the unique neighbor of $y$ in $G$ and let $H^\sigma=G^\sigma-x-y$. Then $y$ is also a pendant vertex of $T_G$ adjacent to $x$. By Lemma~\ref{lem:4.2}, we have
\begin{eqnarray}\label{eq:3.20}
\alpha(T_G)=\alpha(T_G-x)=\alpha(T_H)+1.
\end{eqnarray}

If $x\in W_{\mathscr{C}}$, then the graph $\Gamma_G\bigcup d(G)K_1$ can be obtained from $(T_G-x)\bigcup K_1$ by removing some edges. By Lemma~\ref{lem:3.2}(ii), we get
\begin{eqnarray}\label{eq:3.21}
\alpha(\Gamma_G)+d(G)\geqslant\alpha(T_G-x)+1.
\end{eqnarray}
Now from (\ref{eq:3.20})-(\ref{eq:3.21}) we have $\alpha(\Gamma_G)\geqslant\alpha(T_G-x)-d(G)+1=\alpha(T_G)-d(G)+1$, a contradiction to (iii).

Thus $x$ does not lie on any cycle of $G$. Then $y$ is also a pendant vertex of $\Gamma_G$ adjacent to $x$ and $\Gamma_H=\Gamma_G-x-y$. By Lemma~\ref{lem:4.2} we have
\begin{eqnarray}\label{eq:3.21t}
\alpha(\Gamma_G)=\alpha(\Gamma_H)+1.
\end{eqnarray}
As $x$ does not lie on any cycle of $G$, Lemma~\ref{lem:3.9}(i) implies that
\begin{eqnarray}\label{eq:3.22t}
d(G)=d(H).
\end{eqnarray}
Now from condition (iii) and (\ref{eq:3.20}), (\ref{eq:3.21t})-(\ref{eq:3.22t}), we have $\alpha(T_H)=\alpha(\Gamma_H)+d(H).$ Also note that all cycles of $G$ are cycles of $H$, we conclude that $H^\sigma$ satisfies conditions (i)-(iii). By induction hypothesis we have
\begin{eqnarray}\label{eq:3.22}
sr(H^\sigma)+2\alpha(H)=2(n-2)-2d(H).
\end{eqnarray}
Furthermore, it follows from Lemmas~\ref{lem:3.7} and \ref{lem:4.2} that
\begin{eqnarray}\label{eq:3.23}
sr(G^\sigma)=sr(H^\sigma)+2 \hbox{ and } \alpha(G)=\alpha(H)+1.
\end{eqnarray}
By (\ref{eq:3.22t})-(\ref{eq:3.23}) we have $sr(G^\sigma)+2\alpha(G)=2n-2d(G),$ implying that $G^\sigma$ is lower-optimal.

\setlength{\baselineskip}{15pt}

For ``necessity", let $G^\sigma$ be lower-optimal. By Lemma~\ref{lem:4.6}, the oriented cycles (if any) of $G^\sigma$ are pairwise vertex-disjoint, and each oriented cycle of $G^\sigma$ is odd or evenly-oriented. This implies (i) and (ii).

We proceed by induction on the order $n$ of $G$ to prove (iii). The $n=1$ case is trivial. Suppose that (iii) holds for all lower-optimal oriented graph of order smaller than $n$, and suppose $G^\sigma$ is a lower-optimal oriented graph of order $n\geqslant 2.$

If $T_G$ is an empty graph, then $G^\sigma$ is an odd cycle or an evenly-oriented cycle, in which case (iii) follows immediately. Now suppose $T_G$ has at least one edge, then $T_G$ has at least one pendant vertex, say $y$. Similar to before, either $G$ contains $y$ as a pendant vertex, or $G$ contains a pendant cycle.

\noindent {\bf{Case 1.}}\ $G$ has a pendant vertex $y$.

Let $x$ be the neighbor of $y$ in $G$ and $H^\sigma=G^\sigma-x-y.$ By Lemma~\ref{lem:4.5}, $x$ is not on any cycle of $G$ and $H^\sigma$ is also lower-optimal. Applying induction hypothesis to $H^\sigma$ yields
\begin{eqnarray}\label{eq:3.24}
\alpha(T_H)=\alpha(\Gamma_H)+d(H).
\end{eqnarray}
Since $x$ does not lie on any cycle of $G$, Lemma~\ref{lem:3.9}(i) states that
\begin{eqnarray}\label{eq:3.25}
d(G)=d(H).
\end{eqnarray}
Note that $y$ is also a pendant vertex of $T_G$ (resp. $\Gamma_G$) adjacent to $x$ and $T_H=T_G-x-y$ (resp. $\Gamma_H=\Gamma_G-x-y)$, then by Lemma~\ref{lem:4.2} we have
\begin{eqnarray}\label{eq:3.26}
\alpha(T_G)=\alpha(T_H)+1 \hbox{ and } \alpha(\Gamma_G)=\alpha(\Gamma_H)+1.
\end{eqnarray}
From (\ref{eq:3.24})-(\ref{eq:3.26}) we have
$$
\alpha(T_G)=\alpha(\Gamma_G)+d(G),
$$
as desired.

\noindent {\bf{Case 2.}}\ $G$ has a pendant cycle $C_q$.

Let $x$ be the unique vertex of $C_q$ of degree 3 and $H^\sigma=G^\sigma-C_q^\sigma$. By Lemma~\ref{lem:4.4}(iii), $H^\sigma$ is lower-optimal. Applying the induction hypothesis to $H^\sigma$ yields
\begin{eqnarray}\label{eq:3.27}
\alpha(T_H)=\alpha(\Gamma_H)+d(H).
\end{eqnarray}
From Lemma~\ref{lem:4.4}(ii) we have
\begin{eqnarray}\label{eq:3.28}
\alpha(G)=\alpha(H)+\left\lfloor\frac{q}{2}\right\rfloor.
\end{eqnarray}
Note that $\mathscr{C}_G=\mathscr{C}_H\bigcup{C_q}.$ Together with (\ref{eq:3.28}) and Lemma~\ref{lem:4.6}(iii) we have
\begin{eqnarray}
\alpha(T_G)&=&\alpha(H)+\left\lfloor\frac{q}{2}\right\rfloor-\sum_{C\in \mathscr{C}_G}\left\lfloor\frac{|V_C|}{2}\right\rfloor+d(G)\notag\\
&=&\alpha(H)-\sum_{C\in \mathscr{C}_H}\left\lfloor\frac{|V_C|}{2}\right\rfloor+d(G).\label{eq:3.29}
\end{eqnarray}
Since $H^\sigma$ is lower-optimal, Lemma~\ref{lem:4.6}(iii) states that
\begin{eqnarray}\label{eq:3.30}
\alpha(T_H)=\alpha(H)-\sum_{C\in \mathscr{C}_H}\left\lfloor\frac{|V_C|}{2}\right\rfloor+d(H).
\end{eqnarray}
As $C_q$ is a pendant cycle of $G$, we have
\begin{eqnarray}\label{eq:3.31}
d(G)=d(H)+1.
\end{eqnarray}
Combining (\ref{eq:3.29})-(\ref{eq:3.31}) yields
\begin{eqnarray}\label{eq:3.32}
\alpha(T_G)=\alpha(T_H)+1.
\end{eqnarray}
Note that $\Gamma_G\cong\Gamma_H$, then the required equality $\alpha(T_G)=\alpha(\Gamma_G)+d(G)$ follows from (\ref{eq:3.27}) and (\ref{eq:3.31})-(\ref{eq:3.32}).
This completes the proof.\qed

\subsection{Theorems~\ref{theo:2.4}, \ref{theo:2.5} and \ref{theo:2.6}}\setcounter{equation}{0}

The proofs of  Theorems~\ref{theo:2.4}, \ref{theo:2.5} and \ref{theo:2.6} follow almost directly from Theorem~\ref{theo:2.1}, and are rather similar to each other in nature. Here we only provide the proof of Theorem~\ref{theo:2.4} and leave the rest to the readers.

The \textit{join} of two disjoint graphs $G_1$ and $G_2$, denoted by $G_1\vee G_2$, is the graph obtained from $G_1\cup G_2$ by joining each vertex of $G_1$ to each vertex of $G_2$ by an edge. First we recall the following fact.

\begin{lem}[\cite{29}]
\label{lem:5.1}
Let $G$ be an simple connected graph with $n$ vertices and $m$ edges. Then
$$
\frac{1}{2}\left[(2m+n+1)-\sqrt{(2m+n+1)^2-4n^2}\right]\leqslant\alpha(G)\leqslant\sqrt{n(n-1)-2m+\frac{1}{4}}+\frac{1}{2}.
$$
The equality on the right holds if and only if $G\cong K_{n-\alpha(G)}\vee \alpha(G)K_1$.
\end{lem}

\noindent{\bf Proof of Theorem~\ref{theo:2.4}.}\ Note that for a given simple connected graph $G$ with $|V_G|=n$ and $|E_G|=m$, we have $d(G)=m-n+1$. Together with (\ref{eq:1.1}) and Lemma~\ref{lem:4.1}, we have
\begin{eqnarray*}
sr(G^\sigma)+\alpha(G)\geqslant4n-2m-2-\alpha(G)\geqslant4n-2m-\sqrt{n(n-1)-2m+\frac{1}{4}}-\frac{5}{2}
\end{eqnarray*}
as stated in \eqref{eq:1.2}.

Now we prove the sufficient and necessary conditions for equality in (\ref{eq:1.2}).

\setlength{\baselineskip}{15pt}

``Sufficiency:"
First consider the case that $G\cong S_n$. If $n=1$, then (\ref{eq:1.2}) holds trivially. If $n\geqslant2$,  then we have $sr(G^\sigma)=2$ and $\alpha(G)=n-1$.  Together with the fact that $m=n-1$ we have that equality holds in (\ref{eq:1.2}).

Now we consider the case $G\cong C_3$. By Lemma~\ref{lem:3.6} we have $sr(G^\sigma)=2.$ Note that in this case $\alpha(G)=1$ and $m=n=3$. Hence we have equality in (\ref{eq:1.2}).

``Necessity:" Combining Theorem~\ref{theo:2.1} and Lemma~\ref{lem:4.1} we have that the equality in  (\ref{eq:1.2}) holds if and only if $G^\sigma$ is lower-optimal and $G\cong K_{n-\alpha(G)}\vee \alpha(G)K_1$. Note that the cycles (if any) of $G^\sigma$ are pairwise vertex-disjoint. Hence, $n-\alpha(G)=1$ or $n-\alpha(G)=2$ and $\alpha(G)=1$, which implies $G\cong S_n$ or $G\cong C_3$.

This completes the proof.\qed

\section{Concluding remarks}\label{sec:4}

It is well-known that the AutoGraphiX system determines classes of extremal or near-extremal graphs with a variable
neighborhood search heuristic. As part of a larger study \cite{0001}, the AutographiX2 (AGX2)  \cite{0002,0003,0004}
system was used to study the following type of problems. For each pair of graph invariants $i_1(G)$ and $i_2(G)$, eight bounds of
the following form were considered:
\begin{equation}\label{eq:6.1}
{\underline{b}} \le A\cdot i_1(G)\oplus B\cdot i_2(G)\le \overline{b},
\end{equation}
where $\oplus$ denotes one of the operations $+,-,\times, /$, \, $A, B$ are two constants, while $\underline{b}$ and $\overline{b}$ are, respectively, lower and upper bounding functions. In this paper we considered the invariants $i_1(G)=sr(G^\sigma)$ and $i_2(G)=\alpha(G)$ where $G$ is the underlying graph of $G^\sigma.$
Theorem~\ref{theo:2.1} provides sharp lower bound on $sr(G^\sigma)+2\alpha(G)$; whereas Theorems~\ref{theo:2.4}, \ref{theo:2.5} and \ref{theo:2.6} provide sharp lower bounds on $sr(G^\sigma)+\alpha(G), sr(G^\sigma)-\alpha(G)$ and $sr(G^\sigma)/\alpha(G)$.

It is nature to extend this study through examining the following bounds:
\begin{itemize}
\item sharp upper  bounds on $sr(G^\sigma)+2\alpha(G)$;
\item sharp upper  bounds on $sr(G^\sigma)+\alpha(G), sr(G^\sigma)-\alpha(G)$ and $sr(G^\sigma)/\alpha(G)$;
\item sharp upper and lower bounds on $sr(G^\sigma)\cdot \alpha(G)$.
\end{itemize}

\end{document}